\documentclass[12pt,a4paper]{article}
\usepackage{lipsum}
\usepackage{amsfonts}
\usepackage{graphicx}
\usepackage{epstopdf}
\usepackage[dvipsnames]{xcolor}
\usepackage{mathtools}
\usepackage{amssymb}
\usepackage{esint}
\usepackage{bm}
\usepackage{array}
\usepackage{float}
\usepackage{caption}
\usepackage{subcaption}
\usepackage{enumerate}
\usepackage{algpseudocode}
\usepackage{algorithm}
\usepackage{geometry}
\usepackage[hidelinks]{hyperref} 
\hypersetup{colorlinks=True,
            citecolor=blue} 
\usepackage{placeins}
\usepackage{multirow}
\usepackage{authblk} 
\usepackage{todonotes}
\usepackage{soul}
\usepackage[style = apa]{biblatex}
\DeclareLanguageMapping{english}{english-apa}
\let\cite\parencite

\usepackage[subpreambles=true]{standalone} 
\usepackage{tikz}     

\usepackage{setspace}
\newcommand{\grad}[1]{{\bm{\nabla}{#1}}}
\newcommand{\gradS}[0]{{\bm{\nabla_S}}}
\newcommand{\gradV}[0]{{\bm{\nabla_V}}}

\newcommand{\pfrac}[2]{\frac{\partial{#1}}{\partial{#2}}}

\newcommand{\intOmega}[1]{\int_{\Omega} #1 ~ \text{d} \bm{x}}
\newcommand{\intGamma}[1]{\int_{\Gamma} #1 ~ \text{d}S}

\newcommand{\bracfrac}[2]{\left( \frac{#1}{#2} \right)}
\newcommand{\bracked}[1]{\left( #1 \right)}
\newcommand{\jump}[1]{[[#1]]}

\newcommand{\khat}[0]{\hat{\bm{\mathrm{k}}}}

\let\vec\mathbf

\title{A perfectly matched layer for damping vertically propagating waves in the compressible Boussinesq equations}
\author[1]{Timothy C. Andrews}
\author[2]{Kenneth Duru}
\author[3]{David Lee}

\affil[1]{Department of Climate and Space Sciences and Engineering, University of Michigan, Ann Arbor, MI, USA}
\affil[2]{Department of Mathematical Sciences, University of Texas at El Paso, USA}
\affil[3]{Bureau of Meteorology, Melbourne, Australia}

\begin{document}

\maketitle

\begin{abstract}
    This paper introduces a new application of the perfectly matched layer (PML) for mitigating model top wave reflections in geophysical fluid models. Typically, a strong Laplacian or Rayleigh damping sponge layer is used near the upper boundary, but these often need many vertical levels or a high model top to be sufficiently effective. An advantage of the PML is that, at the continuous level, it is free of wave reflection at the onset of the damping layer. This enables the PML to be effective even with a thin damping layer. We derive PMLs for the linear and nonlinear versions of the Boussinesq equations, which are a simplified model for vertical dynamics in the atmosphere. In the nonlinear system, we define a novel PML that damps perturbations from a hydrostatically balanced reference state. We approximate the PML equations using the compatible finite element method for numerical experiments. First, tests with the linear Boussinesq system show that the PML is more effective than a typical sponge layer in absorbing acoustic waves near the model top. Next, tests in the nonlinear system show that i) the PML can damp acoustic waves even when they are under-resolved by the time discretisation, and ii) the PML can avoid the standing wave pattern caused by model top reflection of the orographic gravity waves. We propose that the PML is worth further development and investigation as a sponge layer alternative in dynamical cores for atmospheric modelling.
\end{abstract}

\section{Introduction}
The absorption of vertically propagating waves near the upper boundary is a key component of numerical models of atmospheric fluid flow. As solving the governing equations on a computer necessitates a truncated vertical domain, some form of artificial boundary is required at the model top. Models with a fixed height coordinate often apply a rigid lid condition that enforces zero vertical velocity at the stationary model top, whilst methods that use floating Lagrangian pressure levels \parencite{lin2004vertically} often apply a free-slip boundary where there is no mass flux through the uppermost pressure level. However, both boundary conditions permit reflections at the model top, causing waves that should have exited the truncated domain to instead propagate downwards and pollute the simulated dynamics. \par
To mitigate model top wave reflections, most atmospheric models introduce a \textit{sponge layer}, which damps the waves that enter the uppermost subdomain \parencite{jablonowski_diffusion}. Sponge layers often use Rayleigh damping, which relaxes the targeted fields back to a reference state. Alternatively, a Laplacian diffusion type damping might be used, or a combination of Laplacian and Rayleigh damping. However, typical sponge layers are often insufficient in idealised tests, or require many levels to provide effective damping; examples include the 10 km sponge layer required for the mesoscale, orographically driven test cases of \textcite{andrews_egusphere-2026-2293}, or the 30 sponge levels recommended for the mountain test case of \textcite{schar2002new}. Other approaches for dealing with vertically propagating waves include the modified Rayleigh damping method of \textcite{klemp2008upper}, or the use of radiation or absorbing boundary conditions, e.g. \textcite{klemp1983upper, klemp1976numerical}. However, these types of boundary conditions are often complicated or infeasible to accurately derive for the governing equations, which contain a spectrum of wave types and frequencies. This work proposes and examines a new possibility for absorbing vertically propagating waves for atmospheric flow problems: the \textit{perfectly matched layer} (PML). \par
Perfectly matched layers have been widely used for electromagnetic applications since their introduction in \textcite{berenger1994perfectly,ChewandWeedon1994}. From there, the PML has been applied to many fields of computational modelling, such as aerospace, seismic, photonic and elastodynamical applications; a summary of important contributions for hyperbolic PDEs can be found in  \textcite{duru2022perfectlymatchedlayerpml}. The PML is a layer that damps all waves that enter, irrespective of their wavelength and angle of incidence. The term \textit{perfectly matched} describes that, in the absence of numerical error, there is no reflection of waves at the boundary between the PML and the rest of the domain. By contrast, Rayleigh damping sponge layers are often very sensitive to the damping strength, with a too large damping coefficient inducing reflection at the damping layer boundary \parencite{klemp1976numerical,andrews_egusphere-2026-2293}. Whilst the PML requires the introduction of auxiliary variables and equations which couple to the original equation set, it typically only requires a small region to provide effective damping. \par
It is significant to note that whilst there are instances of PMLs for geophysical fluid applications, these have been in the context of damping waves at the horizontal boundaries of limited area models of the shallow water equations \parencite{abarbanel2003unsplit, navon2004perfectly, darblade1997methodes}. A related method for the shallow water equations is the absorbing boundary condition approach of \textcite{benacchio2013absorbing}. To our knowledge, the use of a PML to absorb waves in the vertical dimension of atmospheric models has not been investigated. The distinction between horizontal and vertical dynamics is a critical feature of atmospheric dynamics. Where the shallow water equations are the common simplified model for the horizontal dynamics, the equivalent in the vertical dimension is the Boussinesq equations, which examine the impact of stratification on buoyancy oscillations and vertically propagating waves. This motivates our use of the Boussinesq equations in a two-dimensional vertical slice to test the PML for model top damping. \par
We will work with a compressible version of the Boussinesq equations to test the PML for damping two types of atmospheric waves: i) acoustic waves, and ii) orographic gravity waves, which are a subset of internal gravity waves that are excited by flow over a varying surface height. Although many other types of waves may be present in a full atmospheric model, acoustic and orographic gravity waves are two key types that can reflect off the model top. Acoustic waves operate on a very fast timescale, whilst orographic gravity waves propagate much more slowly. In this work, we provide a single PML formulation to damp both acoustic and orographic gravity waves. \par
We derive PML equations for both linear and nonlinear versions of the compressible Boussinesq equations. As the PML is essentially a linear mechanism, for the linear equations we simply apply the existing PML theory for linear hyperbolic systems \cite{duru2022perfectlymatchedlayerpml}. Applications of the PML to nonlinear hyperbolic systems have received much less attention in the literature, aside from the notable exception of \textcite{hu2006construction}. Thus, the derivation of the PML for the nonlinear Boussinesq equations is a novel feature of this study. To succeed in the nonlinear system, we formulate the PML to absorb perturbations from a hydrostatically balanced state. \par
Another novelty of this study is the formulation of compatible finite element approximations of the PML systems. The compatible finite element codebase Gusto is used for numerical tests of the PML for the linear and nonlinear Boussinesq equations. Gusto is a toolbox for solving geophysical fluid equations, designed for rapid prototyping and built from the finite element codebase Firedrake \parencite{FiredrakeUserManual}. We will use simple timestepping methods in Gusto to most clearly highlight the impact of adding PML terms to the original equations. It is important to emphasise that whilst we use the compatible finite element method in this study, the PML equations could be approximated by any numerical method of choice, such as finite difference, finite volume, or continuous and discontinuous Galerkin finite element methods. \par
After this introduction, Section \ref{sec:PML} provides a derivation of the PML for first-order hyperbolic systems, before describing our specific PML for the linear and nonlinear Boussinesq equations. Section \ref{sec:numerical_methods} describes the compatible finite element discretisation used for the numerical tests of section \ref{sec:num_experiments}: A linear Boussinesq test for acoustic wave damping, and a nonlinear test with a mountain to examine acoustic and orographic gravity wave damping. Section \ref{sec:conclusions} ends with conclusions and future research directions.

\section{The perfectly matched layer}
\label{sec:PML}
This section presents a short derivation of the PML for general linear hyperbolic PDEs using the complex coordinate stretching technique of \textcite{ChewandWeedon1994}, as outlined in \textcite{duru2022perfectlymatchedlayerpml}. This approach is then applied to linear Boussinesq equations to form a PML for this system. Subsequently, we extend the PML to the nonlinear Boussinesq equations. Since the PML theory is fundamentally linear, the primary novelty of this section lies in the construction of an accurate PML model in the presence of advective nonlinearities. Last, the parameters of the PML are briefly discussed. 

\subsection{Derivation for first-order hyperbolic PDEs}
\label{subsec:PML_derivation}
We first consider the application of the PML to a linear, first-order, hyperbolic system of equations. The form of PML we consider here is the complex frequency shifted (CFS) variant introduced by \textcite{Kuzuoglu96} and analysed by \textcite{duru2022perfectlymatchedlayerpml}. \par
Consider the initial value, linear partial differential equation (PDE) for a prognostic variable set $\vec{U}$,
\begin{equation}
    \pfrac{\vec{U}}{t} + A \pfrac{\vec{U}}{x} + B \pfrac{\vec{U}}{z} + C \vec{U} = \vec{0}, ~\vec{U}(x,z,t=0) = \vec{U}_0,
\end{equation}

\noindent defined on a vertical slice domain with Cartesian coordinates $(x,z)$, with matrices $A,B,$ and $C$ containing constant coefficients. We seek to damp waves propagating vertically, i.e. in the $z$-coordinate. \par
Consider a Laplace transformation of the prognostic variables from physical space to the complex frequency domain, defined by
\begin{equation}
    \vec{\widetilde{U}}(x,z,s) = \int_{t=0}^{\infty} e^{-st} \vec{U}(x,z,t) ~\text{d}t, ~s \in \mathbb{C}, ~\text{Re} \{s\} > 0.
\end{equation}

\noindent The Laplace-transformed PDE system is
\begin{equation}
   s \vec{\widetilde{U}} + A \pfrac{\vec{\widetilde{U}}}{x} + B \pfrac{\vec{\widetilde{U}}}{z} + C \vec{\widetilde{U}} = \vec{0},
\label{eq:laplace_transformed_pde}
\end{equation}

\noindent where we have tacitly assumed homogenous initial conditions.

Next, we introduce a scaled vertical coordinate of $\zeta=\int_0^z S_z(z') ~\text{d}z'$, with the PML complex metric
\begin{equation}
    S_z(z) = \gamma_z(z) \bracked{1 + \frac{\sigma(z)}{s + \alpha}}.
\label{eq:Sz}
\end{equation}

\noindent Here, $\gamma_z(z)>0$ is a real coordinate stretching \cite{DuKr, Duruthesis2012,duru2022perfectlymatchedlayerpml}, $\sigma(z)\ge 0$ is the PML damping function that is only nonzero in the PML damping region near the model top, and $\alpha \geq 0$ is a real-valued constant defining a complex frequency shift which improves stability of the PML \cite{Kuzuoglu96,DuKr, duru2022perfectlymatchedlayerpml,ElasticDG_PML2019}. Using the coordinate transform (\ref{eq:Sz}) leads to
\begin{equation}
    s \vec{\widetilde{U}} + A \pfrac{\vec{\widetilde{U}}}{x} + \frac{1}{S_z}B \pfrac{\vec{\widetilde{U}}}{\zeta} + C \vec{\widetilde{U}} = \vec{0}.
\end{equation}

Now, we introduce a vector of auxiliary PML variables, $\vec{Q}$. The auxiliary PML variables are defined in the complex frequency domain as
\begin{equation}
    \vec{\widetilde{Q}} = \frac{1}{(s + \alpha)S_z} B \pfrac{\vec{\widetilde{U}}}{\zeta}.
\end{equation}

\noindent We then identify a coupled system of equations for the original and PML variables in the complex frequency domain,
\begin{subequations}
    \begin{align}
        s \vec{\widetilde{U}} + A \pfrac{\vec{\widetilde{U}}}{x} + \frac{1}{\gamma_z}B \pfrac{\vec{\widetilde{U}}}{\zeta} - \sigma \vec{\widetilde{Q}} + C \vec{\widetilde{U}} &= \vec{0}, \\
        s \vec{\widetilde{Q}} - \frac{1}{\gamma_z} B \pfrac{\vec{\widetilde{U}}}{\zeta} + (\sigma + \alpha) \vec{\widetilde{Q}} &= \vec{0},
    \end{align}
\end{subequations}

\noindent using the identity
\begin{equation}
    \frac{1}{S_z} = \frac{1}{\gamma_z} - \frac{1}{S_z} \frac{\sigma}{s + \alpha},
\end{equation}

\noindent which follows from the definition of $S_z$ (\ref{eq:Sz}). Performing the inverse Laplace transform yields the PML equations in physical space:
\begin{subequations}
    \begin{align}
        \pfrac{\vec{U}}{t} + A \pfrac{\vec{U}}{x} + \frac{1}{\gamma_z} B \pfrac{\vec{U}}{\zeta} + C \vec{U} - \sigma \vec{Q} &= \vec{0}, ~\vec{U}(x,z,t=0) = \vec{U}_0, \label{eq:PML_eqs_progs} \\
        \pfrac{\vec{Q}}{t} - \frac{1}{\gamma_z} B \pfrac{\vec{U}}{\zeta} + (\sigma + \alpha) \vec{Q} &= \vec{0}, ~\vec{Q}(x,z,t=0) = \vec{0}. \label{eq:PML_eqs_PML}
    \end{align}
\label{eq:gen_PML_eqs}
\end{subequations}

\noindent Observe that the initial condition for all PML variables is zero. \par
We introduce two important variants of the gradient operator for use in the PML equations, which are a scaled gradient of $\gradS$, and a scaled gradient that only acts in the vertical dimension, $\gradV = (\khat\cdot \gradS)\khat$, where $\khat$ denotes the unit vector in the $z$-coordinate. These operators are
\begin{equation}
    \gradS = \left[ \pfrac{}{x},\frac{1}{\gamma_z}\pfrac{}{z} \right]^T, ~\gradV = \left[ 0,\frac{1}{\gamma_z}\pfrac{}{z} \right]^T.
\label{eq:PML_grad}
\end{equation}

\subsection{A PML for the Boussinesq system}
We now use the theory for linear hyperbolic systems to derive a PML for the linear Boussinesq equations, then discuss the modifications required to apply this to the nonlinear Boussinesq equations. The prognostic variables of the Boussinesq system are a velocity vector of $\vec{u} = [u,w]$, for horizontal velocity $u$ and vertical velocity $w=\vec{u} \cdot \vec{\hat{k}}$, a kinematic pressure $p$ (units of m$^2$ s$^{-2}$), and a buoyancy of $b$ (units of m s$^{-2}$). These pressure and buoyancy fields are deviations from a background, hydrostatically balanced state that is used to derive the Boussinesq equations. There are two parameters in this system: the sound speed $c_{\text{s}}$, and the Brunt-V{\"a}is{\"a}l{\"a} frequency $N$, which measures the timescale of the buoyancy oscillations. 

\subsubsection{Linear Boussinesq equations}
We first apply the PML to a linearised version of the compressible Boussinesq equations, which support acoustic and gravity waves. These equations are described in Section 8.2 of \textcite{Durran} and were used by \textcite{melvin2018choice} to investigate the choice of vertical staggering for the compatible finite element dynamical core, GungHo \parencite{melvin2019gungho_cartesian,melvin2024gungho_spherical}. The linear Boussinesq equations are
\begin{subequations}
\begin{align}
    \pfrac{\vec{u}}{t} + \grad p - b \vec{\hat{k}} &= \vec{0}, \label{subeq:acoustic_buoy_u_vec} \\
    \pfrac{p}{t} + c_{\text{s}}^2 (\bm{\nabla} \cdot \vec{u}) &= 0 \label{subeq:acoustic_buoy_p_eq}, \\
    \pfrac{b}{t} + N^2 (\vec{u} \cdot \vec{\hat{k}}) &= 0. \label{subeq:acoustic_buoy_b_eq}
\end{align}
\label{eq:acoustic_buoy}
\end{subequations}

We now derive a PML for the linear Boussinesq equations (\ref{eq:acoustic_buoy}). As $z$ derivatives only appear in the equations for $\vec{u}$ and $p$, we only require the introduction of two PML variables, $\vec{q_u}$ and $q_p$. This results in the PML equations of
\begin{subequations}
\begin{align}
    \pfrac{\vec{u}}{t} + \bm{\nabla_S} p - b \vec{\hat{k}} - \sigma \vec{q_u} &= \vec{0}, \\
    \pfrac{p}{t} + c_{\text{s}}^2 (\bm{\nabla_S} \cdot \vec{u}) - \sigma q_p &= 0, \\
    \pfrac{b}{t} + N^2 (\vec{u} \cdot \vec{\hat{k}}) &= 0, \\
    \pfrac{\vec{q_u}}{t} - \bm{\nabla_V} p + (\sigma + \alpha) \vec{q_u}&= \vec{0}, \label{eq:acoustic_buoy_q_u_eq} \\
    \pfrac{q_p}{t} - c_{\text{s}}^2 (\bm{\nabla_V} \cdot \vec{u}) + (\sigma + \alpha) q_p&= 0 \label{eq:acoustic_buoy_q_p_eq},
\end{align}
\label{eq:acoustic_buoy_PML}
\end{subequations}

\noindent with the introduction of the scaled PML gradient operators of $\gradS$ and $\gradV$, defined in \eqref{eq:PML_grad}. The PML terms of $- \bm{\nabla_V} p$ and $-c_{\text{s}}^2 (\bm{\nabla_V} \cdot \vec{u})$ damp the vertical wave propagation constructed by the pressure gradient and divergence, which are the key drivers of acoustic wave motion. \par

\subsubsection{Nonlinear Boussinesq equations}
We now move to the nonlinear Boussinesq system. To apply the theory for linear hyperbolic PDEs, the PML is applied to perturbations from a hydrostatically balanced base state of the Boussinesq fields (noting that the Boussinesq fields are themselves perturbations from a background hydrostatic state), derived from linearisation of the nonlinear equations. We use the compressible form of the Boussinesq equations as described by \textcite{gibson2019compatible}, which permits acoustic waves that are not present with the standard Boussinesq approximation. \par
The nonlinear Boussinesq equations are
\begin{subequations}
\begin{align}
    \pfrac{\vec{u}}{t} + (\vec{u} \cdot \grad) \vec{u} + \grad p - b \vec{\hat{k}} &= \vec{0}, \\
    \pfrac{p}{t} + (\vec{u} \cdot \grad) p + c_{\text{s}}^2 (\bm{\nabla} \cdot \vec{u}) &= 0, \\
    \pfrac{b}{t} + (\vec{u} \cdot \grad) b &= 0.
\end{align}
\label{eq:bous}
\end{subequations}

\noindent The main difference from the linear Boussinesq equations (\ref{eq:bous}) is the presence of nonlinear advective transport terms for each prognostic variable. In the buoyancy equation, the advective term replaces the constant buoyancy frequency term of $N^2 (\vec{u} \cdot \khat)$ (\ref{subeq:acoustic_buoy_b_eq}). \par
To linearise the nonlinear Boussinesq equations, we decompose the variables into mean and perturbation components, which are denoted by overbars and primes, respectively:
\begin{equation}
    \vec{u}(x,z) = \vec{u}^{'}, ~p(x,z) = \overline{p}(z) + p^{'}, ~b(x,z) = \overline{b}(z) + b^{'}.
\end{equation}

\noindent Note that the velocity is linearised about a state of rest. \par

The linearisation of the nonlinear Boussinesq equations is then
\begin{subequations}
\begin{align}
    \pfrac{\vec{u}^{'}}{t} + \grad p^{'} - b^{'} \vec{\hat{k}} &= \vec{0}, \\
    \pfrac{p^{'}}{t} + (\vec{u}^{'} \cdot \grad) \overline{p} + c_{\text{s}}^2 (\bm{\nabla} \cdot \vec{u}^{'}) &= 0, \\
    \pfrac{b^{'}}{t} + (\vec{u}^{'} \cdot \grad) \overline{b} &= 0.
\end{align}
\label{eq:bous_linearised}
\end{subequations}

 We apply PML damping to the pressure gradient and divergence terms, as with the linear Boussinesq system, which means that a PML buoyancy variable is again not required. Although the advective terms contain $z$ derivatives corresponding to vertical transport, we neglect these in this study for a couple of reasons. First, the application of the PML to advective terms is often challenging and can be prone to instability, e.g. \textcite{abarbanel1999well,diaz2003stabilized,Hagstrom2003,BecacheEtAL2004}. Second, the vertical velocity is typically much smaller than the horizontal velocity, especially near the model top, where no vertical velocity may be prescribed at the upper boundary. Hence, we neglect the vertical transport terms for this proof-of-concept use of the PML, but will revisit this assumption in future work. \par
The nonlinear Boussinesq PML equations are then
\begin{subequations}
\begin{align}
\pfrac{\vec{u}}{t} + (\vec{u} \cdot \bm{\nabla}) \vec{u} + \bm{\nabla_S} p^{'} - b^{'} \vec{\hat{k}} - \sigma \vec{q_u}&= \vec{0}, \label{eq:PML_bous_mom} \\
    \pfrac{p}{t} + (\vec{u} \cdot \bm{\nabla}) p + c_{\text{s}}^2 (\bm{\nabla_S} \cdot \vec{u}) - \sigma q_p&= 0, \\
    \pfrac{b}{t} + (\vec{u} \cdot \bm{\nabla}) b &= 0, \\
    \pfrac{\vec{q_u}}{t} - \bm{\nabla_V} p^{'} + (\sigma + \alpha) \vec{q_u}&= \vec{0}, \\
    \pfrac{q_p}{t} - c_{\text{s}}^2 (\bm{\nabla_V} \cdot \vec{u}) + (\sigma + \alpha) q_p&= 0.
\end{align}
\label{eq:bous_PML}
\end{subequations}

\noindent Note the explicit use of hydrostatic balance of the reference state, $\bm{\nabla} \overline{p} = \overline{b}\vec{\hat{k}}$, to simplify the momentum equation (\ref{eq:PML_bous_mom}). Setting $\sigma(z) = 0$ restores the original nonlinear Boussinesq equations (\ref{eq:bous}). \par

\subsection{PML parameters}
We now describe the key parameters in the PML, which are summarised in Table \ref{tab:PML_pars}. We consider a PML with a sine-squared damping profile to match the most common choice for Rayleigh damping in dynamical cores. Other choices include the cubic profile that is common in PML literature \parencite{duru2022perfectlymatchedlayerpml}, or quadratic or hyperbolic tangent functions that have alternatively been used in dynamical cores \parencite{jablonowski_diffusion}. \par
Our PML damping function is defined as
\begin{equation}\label{eq:damp_coeff}
\begin{split}
&\sigma(z) = 
\begin{dcases}
0, & \text{if} \quad z < z_{\text{damp}},\\
\sigma_0 \sin^2 \bracked{\frac{\pi}{2}\bracfrac{z - z_{\text{damp}}}{\delta}},  & \text{if}  \quad z \geq z_{\text{damp}} ,
\end{dcases} 
\end{split}
\end{equation}

\noindent with $z_{\text{damp}} = H_z - \delta$ the height at which the PML starts, with $H_z$ the model top height and $\delta$ the thickness of the damping layer, and $\sigma_0$ an inverse damping timescale defined as
\begin{equation}
\sigma_0 = \frac{2c_{\max}}{\delta} \ln \bracfrac{1}{\beta},
\label{eq:sigma_0}
\end{equation}

\noindent with $c_{\text{max}}$ the fastest resolvable wave speed in the domain and $\beta$ a tolerance for the domain truncation error. We set $c_{\text{max}}$ equal to the acoustic wave speed, $c_{\text{s}}$. Our default PML parameters (Table \ref{tab:PML_pars}) and a model top at $H_z = 20 ~\text{km}$ lead to a damping strength of $\sigma_0 = 2.42 ~\text{s}^{-1}$. A good choice of $\sigma_0$ balances being large enough to successfully damp incoming waves, but not too large to degrade the numerical accuracy. \par
For the stretched PML vertical coordinate (\ref{eq:Sz}), we use 
\begin{equation}
    \gamma_z(z) = 1 + \gamma_0 \sigma(z),
\label{eq:gamma_z}
\end{equation}

\noindent with $\gamma_0 \geq 0\in \mathbb{R}$ a real grid stretching factor (with units of s) and $\gamma_0 = 0 ~\text{s}$ reducing to an unstretched grid. A non-zero $\gamma_0$ typically improves the damping ability of the PML, but can lead to accuracy reduction if too large \parencite{martin2009unsplit}. We note that both the unstretched ($\gamma_0=0 ~\text{s}$) and stretched ($\gamma_0>0 ~\text{s}$) PMLs use the same vertical coordinate $z$ in the computational mesh. The grid stretching modifies the transformation in the complex frequency domain, but does not require a modification of the physical domain. Hence, varying the grid stretching factor is a simple modification that only scales derivatives, as reflected in the modified PML operators of (\ref{eq:PML_grad}).  \par
The other key PML parameter is the complex frequency shift (CFS) parameter, $\alpha$, which only appears in the equations for the auxiliary PML variables (\ref{eq:PML_eqs_PML}). Increasing $\alpha$ improves the stability of the PML but reduces its absorption properties \parencite{roden2000convolution}. Following \textcite{duru2025stability}, we set $\alpha=0.05 \sigma_0 ~\text{s}^{-1} = 0.121 ~\text{s}^{-1}$. \par
Our numerical tests, using a finite element method, will use elements of height $\Delta z = 500 ~\text{m}$. This means that for a model top at $H_z = 20 ~\text{km}$ and PML thickness of $\delta = 0.1 H_z = 2 ~\text{km}$, there are only four elements in the vertical extent of the PML. As the auxiliary PML variables only couple to the prognostic equations in the PML region, the extra computation with a PML is restricted to only a small subset of the domain. However, for ease of implementation in these initial studies, we define PML variables that span the entire domain. We will term the region of the domain excluding the PML, $z\in [0,18] ~\text{km}$, the undamped domain.

\begin{table}[]
    \centering
    \caption{Summary of PML parameters and their default values.}
    \begin{tabular}{|c|c|c|}
        \hline
        Parameter & Description & Default \\
        \hline
        $\delta$ & Thickness of the PML & 0.1$H_z ~\text{m}$\\ 
        \hline
        $c_{\text{max}}$ & Fastest wave speed & $c_{\text{s}} = 350 ~\text{m} ~\text{s}^{-1}$ \\
        \hline
        $\beta$ & Tolerance for domain truncation error & $10^{-3}$ \\
        \hline
        $\alpha$ & Complex frequency shift parameter & 0.05$\sigma_0 ~\text{s}^{-1}$ \\
        \hline
        $\gamma_0$ & Real grid stretching factor & 0.25 s \\
        \hline
    \end{tabular}
    \label{tab:PML_pars}
\end{table}

\section{Numerical methods}
\label{sec:numerical_methods}
The PML derived for the Boussinesq equations can be approximated by any numerical method of choice, such as finite difference, finite volume or continuous and discontinuous Galerkin finite element methods. For testing purposes, we will discretise the PML using the compatible finite element method. This section overviews the numerical model and weak forms of the PML Boussinesq equations.

\subsection{Spatial discretisation}
We will work with a compatible finite element model containing quadrilateral cells and next-to-lowest-order elements, which allows for second-order accuracy. A Charney-Phillips vertical staggering is applied, where thermodynamic variables are offset from the dry density in the vertical coordinate. We select the Charney-Phillips grid over the Lorenz grid, where the dry density and thermodynamic variables are instead colocated, to avoid Lorenz grid computational modes \parencite{melvin2018choice,arakawa1996vertical}. An advantage of using the compatible finite element method is that it does not require explicit diffusion for stability \parencite{cotter2012mixed,Cotter23}, meaning that the damping layer can be tested in the absence of additional diffusion mechanisms. \par
\begin{figure}[htpb]
    \centering
     \includegraphics[width=\textwidth]{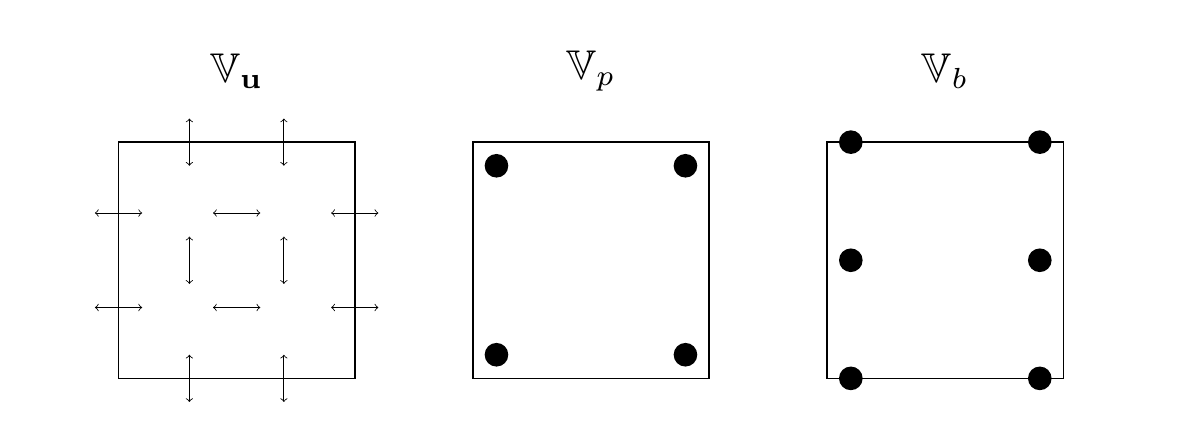}
    \caption[]{Next-to-lowest order finite element function spaces used in this work. The vector Raviart-Thomas function space of $\mathbb{V}_{\vec{u}}$ enforces continuity of the normal component of the vector field, as shown by the direction of the arrows at the cell facets. The two scalar function spaces of $\mathbb{V}_{p}$ and $\mathbb{V}_{b}$, which are used with the Charney-Phillips vertical staggering, have degrees of freedom denoted by circles, with circles at the facets for $\mathbb{V}_{b}$ implying $\mathcal{C}^{0}$ continuity between adjacent vertical elements.}
    \label{fig:function_spaces}
\end{figure}

The quadrilateral elements are constructed as the tensor product, $\otimes$, of a one-dimensional finite element in the horizontal with a one-dimensional finite element in the vertical \parencite{mcrae2016automated}. For degree $l$, we denote a discontinuous element as $\text{d}Q_l$ and a continuous element as $Q_l$. We then use the following types of function space (visualised in Figure \ref{fig:function_spaces}):
\begin{itemize}
    \item $\mathbb{V}_{\vec{u}}:= \text{RT}_1$, the next-to-lowest-order Raviart-Thomas function space on quadrilateral elements, which enforces continuity of normal components of the vector field at the facets, whilst the tangential component is permitted to be discontinuous. The horizontal component of the vector field is represented by continuous order 2 polynomials in the horizontal and discontinuous order 1 polynomials in the vertical ($Q_2 \otimes \text{d}Q_1$). Conversely, the vertical component is represented by discontinuous order 1 polynomials in the horizontal and continuous order 2 polynomials in the vertical ($\text{d}Q_1 \otimes Q_2$). 
    \item $\mathbb{V}_{p} := \text{d}Q_1 \otimes \text{d}Q_1$, the fully discontinuous next-to-lowest order function space. 
    \item $\mathbb{V}_{b} := \text{d}Q_1 \otimes Q_2$, the next-to-lowest-order Charney-Phillips function space for thermodynamic variables, i.e. buoyancy, which is discontinuous in the horizontal and continuous in the vertical.
\end{itemize}

\noindent For a Charney-Phillips staggering, we take $(\vec{u},p,b) \in (\mathbb{V}_{\vec{u}}, \mathbb{V}_{p}, \mathbb{V}_{b})$. We define the PML variables in the same function space as their prognostic counterpart, so $\vec{q_u} \in \mathbb{V}_{\vec{u}}, q_p \in \mathbb{V}_{p}$. \par

\subsection{Weak forms}
Here we present the weak forms used for our compatible finite element discretisations. We will consider first the weak form of the PML for the linear Boussinesq equations, then proceed to the nonlinear Boussinesq equations.
\subsubsection{Linear Boussinesq equations}
To obtain a weak form of the linear Boussinesq PML equations (\ref{eq:bous_PML}), we take the inner product of each equation with a corresponding test function, using $\tau$ to denote prognostic test functions and $\chi$ for PML test functions. Integrating over the domain, $\Omega$, arrives at
\begin{subequations}
\begin{align}
    \intOmega{\bm{\tau_u} \cdot \pfrac{\vec{u}}{t}} - \intOmega{p  \grad \cdot (\bm{\tau_u} \odot \bm{\gamma_S})} - \intOmega{\bm{\tau_u} \cdot b \hat{\vec{k}}}& \nonumber \\
    - \intOmega{\bm{\tau_u} \cdot \sigma \vec{q_u}}&= 0, ~\forall \bm{\tau_u} \in \mathbb{V}_{\vec{u}}, \label{eq:lin_bous_mom_weak}  \\
    \intOmega{\tau_p \pfrac{p}{t} } + \intOmega{c_{\text{s}}^2 \tau_p (\bm{\nabla_S} \cdot \vec{u})} - \intOmega{\sigma \tau_p q_p } &= 0, ~\forall \tau_p \in \mathbb{V}_{p}, \label{eq:acoustic_buoy_weak_p_eq} \\
    \intOmega{\tau_b \pfrac{b}{t}} + \intOmega{\tau_b N^2 (\vec{u} \cdot \hat{\vec{k}})} &= 0, ~\forall \tau_b \in \mathbb{V}_{b}, \label{eq:acoustic_buoy_weak_b_eq} \\
    \intOmega{\bm{\chi_u} \cdot \pfrac{\vec{q_u}}{t}} + \intOmega{p \grad \cdot (\bm{\chi_u} \odot \bm{\gamma_V)}} + \intOmega{\bm{\chi_u} \cdot (\sigma + \alpha) \vec{q_u}} &= 0, ~\forall \bm{\chi_u} \in \mathbb{V}_{\vec{u}}, \\
    \intOmega{ \chi_p \pfrac{q_p}{t}  } - \intOmega{c_{\text{s}}^2 \chi_p (\bm{\nabla_V} \cdot \vec{u})} + \intOmega{(\sigma + \alpha) \chi_p q_p } &= 0, ~\forall \chi_p \in \mathbb{V}_{p},
\end{align}
\label{eq:bous_PML_weak}
\end{subequations}

\noindent where $\bm{x}$ denotes the spatial coordinate. The application of integration by parts on the pressure gradient means that the PML grid stretching term $\gamma_z(z)$ is operated upon by the divergence operator. We express this by introducing PML scaling terms of $\bm{\gamma_S} = [1, \gamma_z]^T$ and $\bm{\gamma_V} = [0, \gamma_z]^T$, and the Hadamard (element-wise) product of two vectors, denoted as $\odot$ (note that $\gradS = \bm{\nabla} \odot \bm{\gamma_S}, \gradV = \bm{\nabla} \odot \bm{\gamma_V}$). There are no surface terms left after integration by parts due to the continuity of normal components of vector fields in the Raviart-Thomas function space (Figure \ref{fig:function_spaces}). \par
The linear Boussinesq equations conserve the energy
\begin{equation}
    \mathcal{H} = \frac{1}{2} \intOmega{\bracked{\vec{u} \cdot \vec{u} + c_{\text{s}}^{-2} p^2 + N^{-2} b^2}},
\label{eq:energy}
\end{equation}

\noindent which can be decomposed into kinetic energy ($E_{\text{K}}$), internal energy ($E_{\text{I}}$) and potential energy ($E_{\text{P}}$) components:
\begin{subequations}
    \begin{align}
        E_{\text{K}} &= \frac{1}{2} \intOmega{\vec{u} \cdot \vec{u}}, \label{eq:E_K} \\
        E_{\text{I}} &= \frac{1}{2} \intOmega{c_{\text{s}}^{-2} p^2 }, \label{eq:E_I} \\
        E_{\text{P}} &= \frac{1}{2} \intOmega{N^{-2} b^2}. \label{eq:E_P}
    \end{align}
\label{eq:bous_energies}
\end{subequations}

To examine energy conservation in the Boussinesq equations, we take the weak form PML system and scale the pressure equation (\ref{eq:acoustic_buoy_weak_p_eq}) by $c_{\text{s}}^{-2}$ and the buoyancy equation by $N^{-2}$ (\ref{eq:acoustic_buoy_weak_b_eq}), so that
\begin{subequations}
\begin{align}
    \intOmega{\bm{\tau_u} \cdot \pfrac{\vec{u}}{t}} - \intOmega{p  \grad \cdot (\bm{\tau_u} \odot \bm{\gamma_S})} 
    - \intOmega{\bm{\tau_u} \cdot b \hat{\vec{k}}} &= \intOmega{\bm{\tau_u} \cdot \sigma\vec{q_u}}, &~\forall \bm{\tau_u} \in \mathbb{V}_{\vec{u}}, \\
    \intOmega{c_{\text{s}}^{-2} \tau_p \pfrac{p}{t} } + \intOmega{ \tau_p (\bm{\nabla_S} \cdot \vec{u})} &=  \intOmega{ c_{\text{s}}^{-2}\sigma \tau_p q_p }, &~\forall \tau_p \in \mathbb{V}_{p}, \\
    \intOmega{N^{-2} \tau_b \pfrac{b}{t}} + \intOmega{\tau_b (\vec{u} \cdot \hat{\vec{k}})} &= 0, &~\forall \tau_b \in \mathbb{V}_{b}.
\end{align}
\end{subequations}

\noindent Choosing test functions of $\bm{\tau_u} = \vec{u}, ~\tau_p=p,$ and $\tau_b=b$ shows that the energy (\ref{eq:energy}) is conserved in time in the case of no PML, i.e. $d\mathcal{H}/dt =0$ when $\sigma=0$. When including the PML, we have a source term which modifies the energy as
\begin{equation}
    \frac{\mathrm{d} \mathcal{H}}{\mathrm{d} t} = P_{\text{PML}} = \intOmega{\sigma (\vec{u} \cdot \vec{q_u} + c_{\text{s}}^{-2} p q_p)},
\label{eq:P_PML}
\end{equation}

\noindent with $P_{\text{PML}}$ denoting the PML power. Energy-conserving formulations for mixed finite element discretisations have been presented previously for other geophysical systems such as the compressible Euler equations \parencite{Lee21,LeePalha21,Cotter23}.

\subsubsection{Nonlinear Boussinesq equations}
For the weak form of the nonlinear Boussinesq PML equations (\ref{eq:bous_PML}), discontinuous Galerkin (DG) upwinding is applied to the advective terms \parencite{cockburn2001runge,bendall2023improving},
\begin{subequations}
\begin{align}
    \intOmega{\bm{\tau_u} \cdot \pfrac{\vec{u}}{t}} + \intGamma{ (\vec{u} \cdot \vec{\hat{n}^+}) \vec{u}^{\dagger} \cdot [[\bm{\tau_u}]]} & \nonumber \\
    - \intOmega{\vec{u} \cdot [\nabla \cdot (\bm{\tau_u} \otimes \vec{u})]}
     - \intOmega{p^{'}  \bm{\nabla} \cdot (\bm{\tau_u} \odot \bm{\gamma_S})} \nonumber \\
     - \intOmega{ \bm{\tau_u} \cdot b^{'}\khat} - \intOmega{ \bm{\tau_u} \cdot \sigma\vec{q_u}} &= 0, ~\forall \bm{\tau_u} \in \mathbb{V}_{\vec{u}}, \label{eq:nl_bous_mom_weak}\\
    \intOmega{\tau_p \pfrac{p}{t} } + \intGamma{(\vec{u}  \cdot \vec{\hat{n}}^+) p^{\dagger} [[\tau_p]]} - \intOmega{p \bm{\nabla} \cdot (\tau_p \vec{u}) } & \nonumber \\
     + \intOmega{c_{\text{s}}^2 \tau_p (\bm{\nabla_S} \cdot \vec{u})} - \intOmega{\sigma \tau_p q_p } &= 0, ~\forall \tau_p \in \mathbb{V}_{p}, \\
    \intOmega{\tau_b \pfrac{b}{t} } + \intGamma{(\vec{u} \cdot \vec{\hat{n}}^+) b^{\dagger} [[\tau_b]]}
    - \intOmega{b \bm{\nabla} \cdot (\tau_b \vec{u}) } &= 0, ~\forall \tau_b \in \mathbb{V}_{b}, \\
    \intOmega{\bm{\chi_u} \cdot \pfrac{\vec{q_u}}{t}} + \intOmega{p^{'} \grad \cdot (\bm{\chi_u} \odot \bm{\gamma_V)}} + \intOmega{(\sigma + \alpha) \bm{\chi_u} \cdot \vec{q_u}} &= 0, ~\forall \bm{\chi_u} \in \mathbb{V}_{\vec{u}}, \\
    \intOmega{ \chi_p \pfrac{q_p}{t}  } - \intOmega{c_{\text{s}}^2 \chi_p (\gradV \cdot \vec{u}) } + \intOmega{(\sigma + \alpha) \chi_p q_p } &= 0, ~\forall \chi_p \in \mathbb{V}_{p},
\end{align}
\label{eq:bous_PML_vector_weak}
\end{subequations}

\noindent with $\Gamma$ denoting the set of interior facets, and $^\dagger$ denoting the upwind value on a facet, e.g.
\begin{equation} \label{def:upwind}
    \vec{q}^\dagger := \left\lbrace 
    \begin{matrix}
    \vec{q}^+ & \mathrm{if} \ \vec{u} \cdot \hat{\vec{n}}^+ \geq 0, \\
    \vec{q}^- & \mathrm{if} \ \vec{u} \cdot \hat{\vec{n}}^+ < 0,
    \end{matrix}\right.
    \end{equation} 

\noindent where the `+' and `-' superscripts refer to arbitrarily labelled sides of the facet, $\hat{\vec{n}}$ are outwardly facing normals, and $\vec{q}$ could be a scalar or vector. The double square brackets denote the jump of a field across a facet,
    \begin{equation}
        \jump{\vec{q}} = \vec{q}^+ - \vec{q}^-.
    \end{equation}

\subsection{Rayleigh damping sponge layer}
\label{subsec:sponge}
To evaluate the PML for damping vertically propagating waves, we compare its performance to that of a Rayleigh damping sponge layer. Following \textcite{cotter2023compatible}, we introduce a sponge term of $\mu \khat (\vec{u} \cdot \khat)$. The corresponding weak form of 
\begin{equation}
    \intOmega{\mu(\bm{\tau_u} \cdot \khat)(\vec{u} \cdot \khat)},
\label{eq:sponge_term}
\end{equation}

\noindent is added to the left-hand sides of the momentum equations, so (\ref{eq:lin_bous_mom_weak}) for the linear Boussinesq system and (\ref{eq:nl_bous_mom_weak}) for the nonlinear Boussinesq system. This Rayleigh damping term acts to damp the vertical velocity back to zero. Rayleigh damping can also be applied to horizontal velocities \parencite{harris2021scientific_FV3_doc,polvani2002tropospheric,klemp1976numerical,jablonowski_diffusion}, with damping of the vertical velocity being more common in nonhydrostatic models, such as MPAS \parencite{skamarock2012multiscale}, ICON \parencite{zangl2015icon}, and GungHo \parencite{melvin2024gungho_spherical}. \par
The sponge layer damping function $\mu(z)$ is defined analogously to $\sigma(z)$ for the PML (\ref{eq:damp_coeff}), with an inverse damping timescale of $\mu_0$ used in place of $\sigma_0$. For the acoustic tests, we examine sponge strengths for values of $\mu_0$ separated by $\Delta \mu_0 = 0.5 ~\text{s}^{-1}$. We note that more vertical levels can improve the efficacy of a sponge layer, but we use the same thickness as the PML for fair comparison, which is four elements in the vertical. The use of only a few models levels in the sponge layer is standard for dynamical cores, where it is important to minimise the additional computation, e.g. typically three layers are used in the finite volume cubed-sphere (FV3) model \parencite{harris2021scientific_FV3_doc} and the Community Atmosphere Model (CAM) \parencite{neale2010description} version of the spectral element model \parencite{dennis2012cam,lauritzen2018ncar}. \par

\section{Numerical experiments}
\label{sec:num_experiments}
We now describe tests of the PML in the linear and nonlinear Boussinesq equations, to show the effective PML absorption properties for outgoing waves and verify numerical accuracy and stability. In particular, we compare the performance of the PML to the Rayleigh damping sponge layer that is the default for dynamical cores. \par
In these tests, we use a sound speed of $c_{\text{s}} = 350 ~\text{m} ~\text{s}^{-1}$. The Brunt-V{\"a}is{\"a}l{\"a} frequency is set to a constant $N=0.01 ~\text{s}^{-1}$, which mimics an isothermal atmosphere. The domain is $[0,L_x] \times [0,H_z]$, with length $L_x$ = 100 km and height $H_z$ = 20 km. Next-to-lowest order elements of width $\Delta x = 500 ~\text{m}$ and height $\Delta z = 500 ~\text{m}$ are used. 

\subsection{Linear Boussinesq equations}
We begin with the linear system, where we test the damping of vertically propagating acoustic waves.

\subsubsection{Test setup}
We consider an initial state that contains a perturbation to the pressure field. Often, the buoyancy field is instead perturbed in Boussinesq systems, most famously in the gravity wave test of \textcite{skamarock1994efficiency}. However, the buoyancy perturbation excites predominantly horizontally propagating gravity waves that are much less prone to model top reflection in this configuration. \par
The pressure perturbation takes the form of a compactly supported Gaussian function at the centre of the domain,
\begin{equation}
p(x,z,t=0)=
    \begin{cases}
        p_0 \exp\bracked{-\bracked{\frac{(x-L_x/2)^2+(z-H_z/2)^2}{d^2}}}, &\text{if} ~x \in [0.3L_x, 0.7L_x] ~\text{and} ~z \in [0.3H_z, 0.7H_z], \\
        0, &\text{else}.
    \end{cases}
\end{equation}

\noindent The maximum pressure perturbation magnitude is $p_0 = 10 ~\text{Pa}$ and the Gaussian half-width is $d = 1 ~\text{km}$. The initial velocity and buoyancy fields are zero. A flat surface at $z_{\text{s}}=0$ is used. The acoustic waves take 28.6 s to travel $10 ~\text{km}$, i.e. from the centre of the domain to the model top or surface. An explicit RK4 timestepper is used with $\Delta t = 0.25 ~\text{s}$, which leads to an acoustic wave Courant number of $\textrm{Cr}_{\text{ac}}= c_{\text{s}} \Delta t/\Delta x=0.175$. The simulation length is $T_{\text{end}} = 100 ~\text{s}$. \par
\begin{figure}[htpb]
    \centering
    \includegraphics[width=0.75\linewidth]{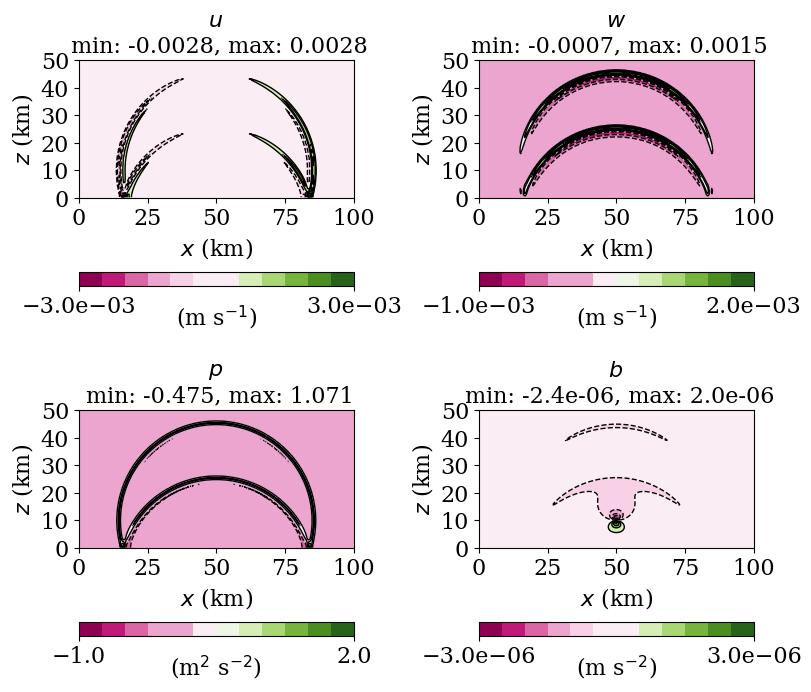}
    \caption{A reference solution at $t=100$ s of the linear Boussinesq test with a higher model top of 50 km.}
    \label{fig:bous_ref}
\end{figure}

We compute a reference solution for quantitative comparisons of no model top damping, a sponge layer, or a PML. The reference solution applies no model top damping, but uses a higher domain with $H_z = 50 ~\text{km}$ to avoid wave reflections. Figure \ref{fig:bous_ref} shows the end state of the reference solution. There are two wave fronts, as the initial pressure perturbation excites upwardly and downwardly travelling components. The downwardly travelling component reflects off the bottom surface and then travels upwards, thus ending at a lower height than the upwardly travelling crest. In the next section, errors relative to this reference solution are computed as the total $L_2$ error of the vertical velocity within the undamped domain of $z \in [0,18] ~\text{km}$.

\subsubsection{Results}
\begin{figure}[htpb]
    \centering
    \includegraphics[width=0.95\linewidth]{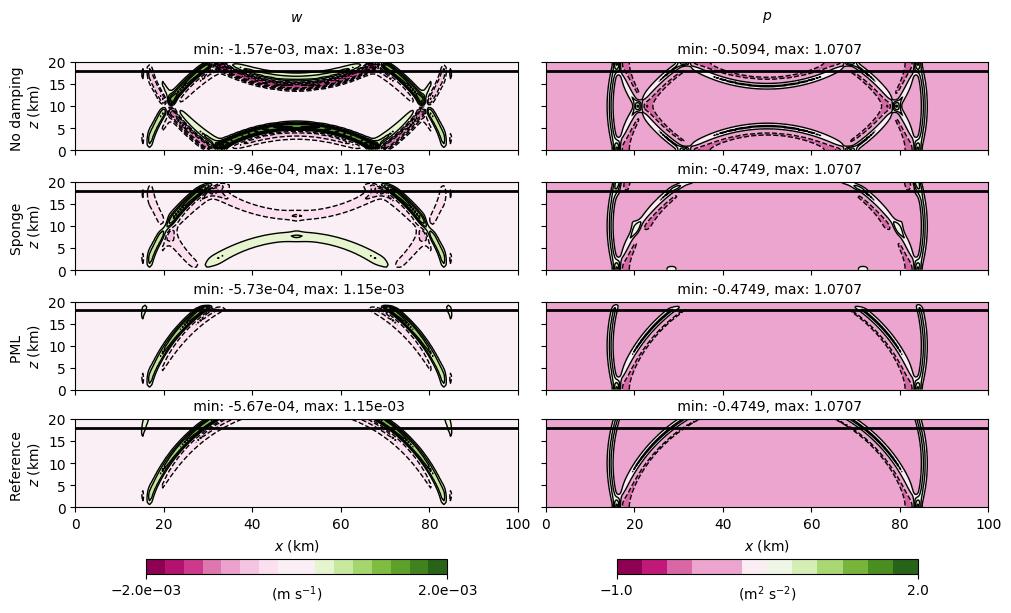}
    \caption{Vertical velocity (left column) and pressure (right column) fields from the linear Boussinesq test after 100 s. Solutions are shown with no model top damping (top row), a sponge layer (second row), a PML with $\gamma_0 = 0.25$ s (third row), and from the high model top reference (bottom row). The black horizontal line indicates the separation between the undamped domain of $z \in [0,18] ~\text{km}$ and the damping layer of $z \in [18, 20] ~\text{km}$ (if one is used). The field extrema are computed only in the undamped domain.}
    \label{fig:acoustic_buoyancy_compare_fields}
\end{figure}

The vertical velocity and pressure fields from simulations with no damping, a sponge layer, and the PML with $\gamma_0=0.25 ~\text{s}$ are compared with the reference solution in Figure \ref{fig:acoustic_buoyancy_compare_fields}. The best sponge strength leading to the smallest vertical velocity error is identified as $\mu_0 = 1.5 ~\text{s}^{-1}$, which is weaker than the PML coefficient of $\sigma_0 = 2.42 ~\text{s}^{-1}$ (\ref{eq:sigma_0}). In Figure \ref{fig:acoustic_buoyancy_compare_fields}, the undamped solution has clear model top reflections from both the upwardly and downwardly propagating parts of the acoustic wave, the latter after it has reflected off the bottom surface. With a sponge layer, the reflection is weaker and the amplitudes of $w$ are reduced, but there is still a noticeable spurious impact on the dynamics. With the PML, both the upwardly and downwardly propagating acoustic crests are successfully damped once they enter the PML, leaving no visual evidence of reflection. Accordingly, the solution is visually indistinguishable from the reference solution in the undamped domain. \par

\begin{figure}
    \centering
    \includegraphics[width=0.7\linewidth]{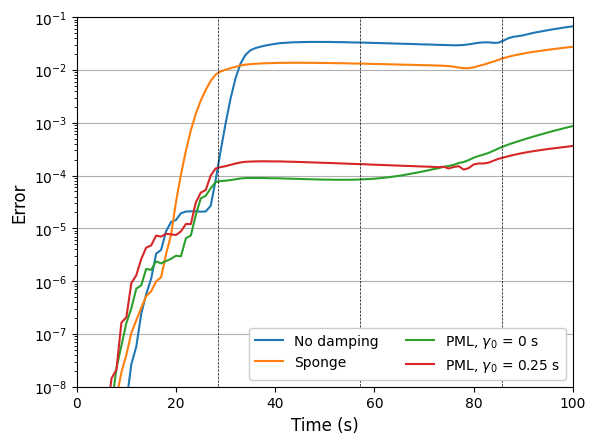}
    \caption{$L_2$ vertical velocity error in the undamped domain, $z \in [0,18] ~\text{km}$, for the linear Boussinesq test in simulations with no model top damping, a sponge layer, and a PML (with and without grid stretching). The vertical dashed lines indicate time intervals of $28.6 ~\text{s}$, which is the time taken for the acoustic wave to travel 10 km.}
    \label{fig:acoustic_buoyancy_errors}
\end{figure}

The improvement with the PML is quantified through a time series of the vertical velocity error in Figure \ref{fig:acoustic_buoyancy_errors}. Both the unstretched ($\gamma_0=0 ~\text{s}$) and stretched ($\gamma_0 \neq 0~\text{s}$) PML solutions are two orders of magnitude more accurate than solutions with no damping or a sponge layer. The use of a grid stretching factor of $\gamma_0=0.25 ~\text{s}$ further reduces the error compared to the unstretched PML when $t > 75 ~\text{s}$; this stretching factor was identified as the choice with lowest vertical velocity error from $\gamma_0 \in \{0, 0.05, 0.1, 0.25, 0.5 \} ~\text{s}$. Hence, the use of real grid stretching can improve the damping ability of the PML, although making this too large begins to degrade the accuracy. \par

\begin{figure}[htpb]
    \centering
    \includegraphics[width=0.9\linewidth]{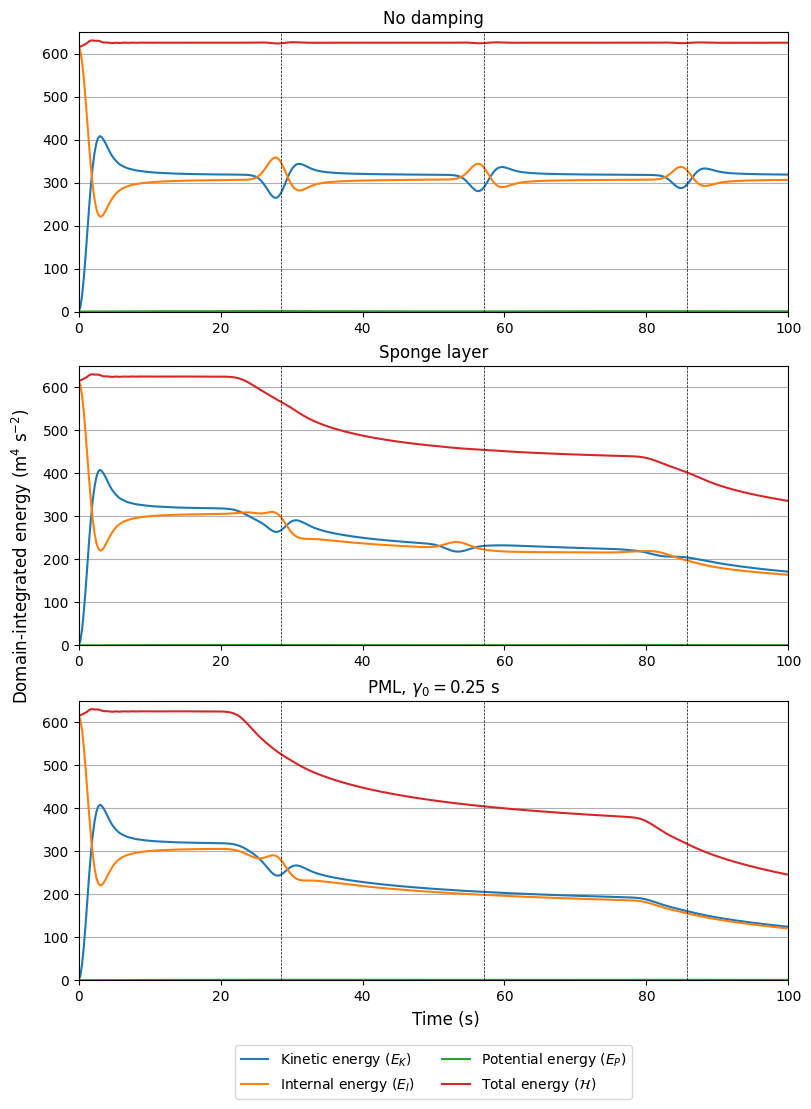}
    \caption{Time series of energy components in the linear Boussinesq test, with no model top damping, a sponge layer, and a stretched PML. The vertical dashed lines indicate time intervals of $28.6 ~\text{s}$, which is the time taken for the acoustic wave to travel 10 km.}
    \label{fig:acoustic_buoyancy_energies}
\end{figure}

We now examine the energetics of the linear Boussinesq system (\ref{eq:energy}) in Figure \ref{fig:acoustic_buoyancy_energies},,. The initial condition only contains internal energy (due to the linear equation set solving for perturbations around a stably stratified state) which is primarily exchanged with kinetic energy during the simulation. In the undamped case, there are three bumps in the kinetic and internal energies after the initial disturbance, which are due to i) the upwardly travelling crest of the acoustic wave hitting the model top and the downwardly travelling crest hitting the bottom surface ($t \approx 28.6$ s), ii) the upwardly and downwardly travelling crests meeting in the middle of the domain ($t \approx 57.1$ s), and iii), the initially downwardly travelling wave hitting the model top ($t \approx 85.7$ s). \par
The undamped solution maintains a fairly constant total energy throughout the simulation, which is a consequence of a conserved energy in the linear Boussinesq equations (\ref{eq:energy}). With a sponge layer, there is a reduction in total energy once the upwardly travelling crest enters the damping region. The second and third bumps in $E_K$ and $E_I$ are still present with a sponge, although with a smaller amplitude than in the undamped case. With the PML, there is again a reduction in total energy as the upwardly travelling crest enters the damping region. A faster energy drop off with the PML compared to the sponge indicates a more effective removal of acoustic wave energy. The first bump in the energies is present, as there is still a reflection of the downwardly travelling crest at the bottom surface, but there are no second and third bumps, indicating the absence of wave reflection at the model top. \par

\begin{figure}[htpb]
    \centering
    \includegraphics[width=0.8\linewidth]{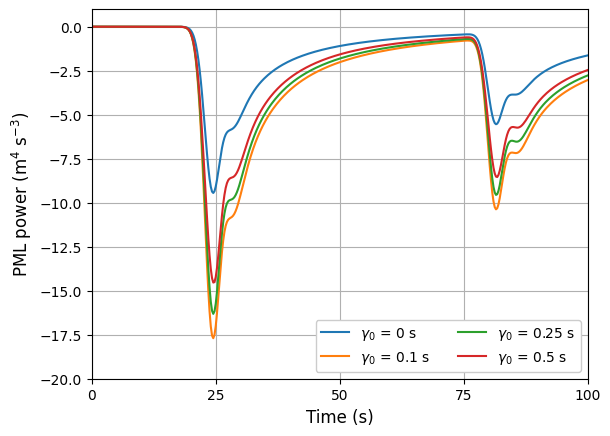}
    \caption{Time series of the PML power term, $P_{\text{PML}}$ (\ref{eq:P_PML}), in the linear Boussinesq test. }
    \label{fig:acoustic_buoyancy_PML_power}
\end{figure}

Lastly, Figure \ref{fig:acoustic_buoyancy_PML_power} shows the rate of energy removal in the PML for different grid stretching factors. The PML power (\ref{eq:P_PML}) is zero until the upwardly travelling crest enters the PML region and is damped. The PML power peaks when damping this first acoustic wave. There is also a second region of strong PML power when the downwardly travelling crest reaches the PML region after its reflection off the bottom surface. The most energy removed is with $\gamma_0=0.1 ~\text{s}$, although $\gamma_0 = 0.25 ~\text{s}$ leads to the smallest vertical velocity error in this test.

\FloatBarrier

\subsection{Nonlinear Boussinesq equations}
\label{sec:Boussinesq}
We now move to the nonlinear Boussinesq equations, where we test the damping of both acoustic and orographic gravity waves. 

\subsubsection{Test setup}
We apply an analytical initial condition of
\begin{subequations}
    \begin{align}
        \vec{u}(x,z,t=0)=\overline{\vec{u}} &= [10,0]^T, \\
        b(x,z,t=0)=\overline{b}(z) &= N^2 z, \\
        p(x,z,t=0)=\overline{p}(z) &= \frac{N^2 z^2}{2} - p_0,
    \end{align}
\end{subequations}

\noindent with $p_0 = 200 ~\text{hPa}$. This initial condition satisfies the Boussinesq hydrostatic balance of $\grad p=b\vec{\hat{k}}$. Alternatively, one could enforce hydrostatic balance numerically with the method of \textcite{natale2016compatible}, which iteratively computes a $\overline{p}$ that is in discrete hydrostatic balance with $\overline{b}$ in each column to some error tolerance. We choose an analytical expression as this enables the computation of a reference solution with a higher model top, whereas the discrete hydrostatic balance method leads to a different $\overline{p}(z)$ with higher columns. \par
We use a varied surface height to enable orographic gravity waves in this test. We define a surface height containing a Gaussian mountain,
\begin{equation}
    z_{\text{s}} = h_0 \exp \bracked{-\bracfrac{x-x_{\text{c}}}{a}^2 },
\end{equation}

\noindent with a mountain height of $h_0 = 1 ~\text{km}$, and Gaussian half-width of $a = 10 ~\text{km}$. The mountain is centred at the middle of the domain, $x_{\text{c}}=L_x/2$. The mountain is sufficiently high to enter the nonlinear `flow-around' regime, which leads to more complicated dynamics than the linear `flow-over' regime of small mountains \parencite{smith1989hydrostatic,lin2007mesoscale}.  \par 
We use a hybrid terrain-following vertical coordinate of
\begin{equation}
    z = \overline{z} + \frac{H_z-\overline{z}}{H_z} z_{\text{s}},
\label{eq:hybrid_z_coord}
\end{equation}

\noindent with $\overline{z}$ a reference vertical coordinate with constant height levels evenly spaced by $\Delta \overline{z} = 500 ~\text{m}$. This choice ensures that $z=z_{\text{s}}$ at the surface and the model top at $z=H_z$ is flat. \par
The same initial condition is used for the acoustic and orographic gravity wave configurations of this test. The acoustic wave is generated through an incompatibility of the horizontal flow in the initial condition with the impermeability boundary condition of $ \vec{u} \cdot \vec{\hat{n}}=0$ at the bottom surface. Unlike the linear Boussinesq test, the acoustic wave starts at the surface and only contains an upwardly travelling crest. The orographic gravity wave is generated from the horizontal flow over the mountain and propagates on a much slower timescale than the acoustic wave. To specifically test the damping of each type of wave, we vary the time discretisation and simulation length:
\begin{itemize}
    \item For the acoustic test: Explicit RK4 timestepping, $\Delta t = 0.25 ~\text{s}, ~T_{\text{end}} = 150 ~\text{s}$. The acoustic wave Courant number is $\textrm{Cr}_{\text{ac}}=0.175$.
    \item For the orographic gravity wave test: Implicit midpoint timestepping (also called the Crank-Nicolson scheme or trapezoidal rule), $\Delta t = 10 ~\text{s}, ~T_{\text{end}} = 10,000 ~\text{s}$. The acoustic Courant number is $\textrm{Cr}_{\text{ac}}=7$ and advective Courant number is $\textrm{Cr}_{\text{adv}}=u_0 \Delta t/\Delta x = 0.2$, with $u_0 = 10 ~\text{m} ~\text{s}^{-1}$. 
\end{itemize}

\subsubsection{Acoustic test results}
For the acoustic test, we again compute a reference solution with a higher model top of 50 km and no damping layer. The reference solution uses the terrain following vertical coordinate (\ref{eq:hybrid_z_coord}) for $z \in [0, 20] ~\text{km}$ with equally spaced levels of element height $\Delta z=500 ~\text{m}$ for $z\in [20,50] ~\text{km}$. With the simulation time of $T_{\text{end}} = 150 ~\text{s}$, the acoustic waves travel 52.5 km and do not interfere with the solution in the undamped domain of $z \in [0, 18] ~\text{km}$. \par
\begin{figure}
    \centering
    \includegraphics[width=\linewidth]{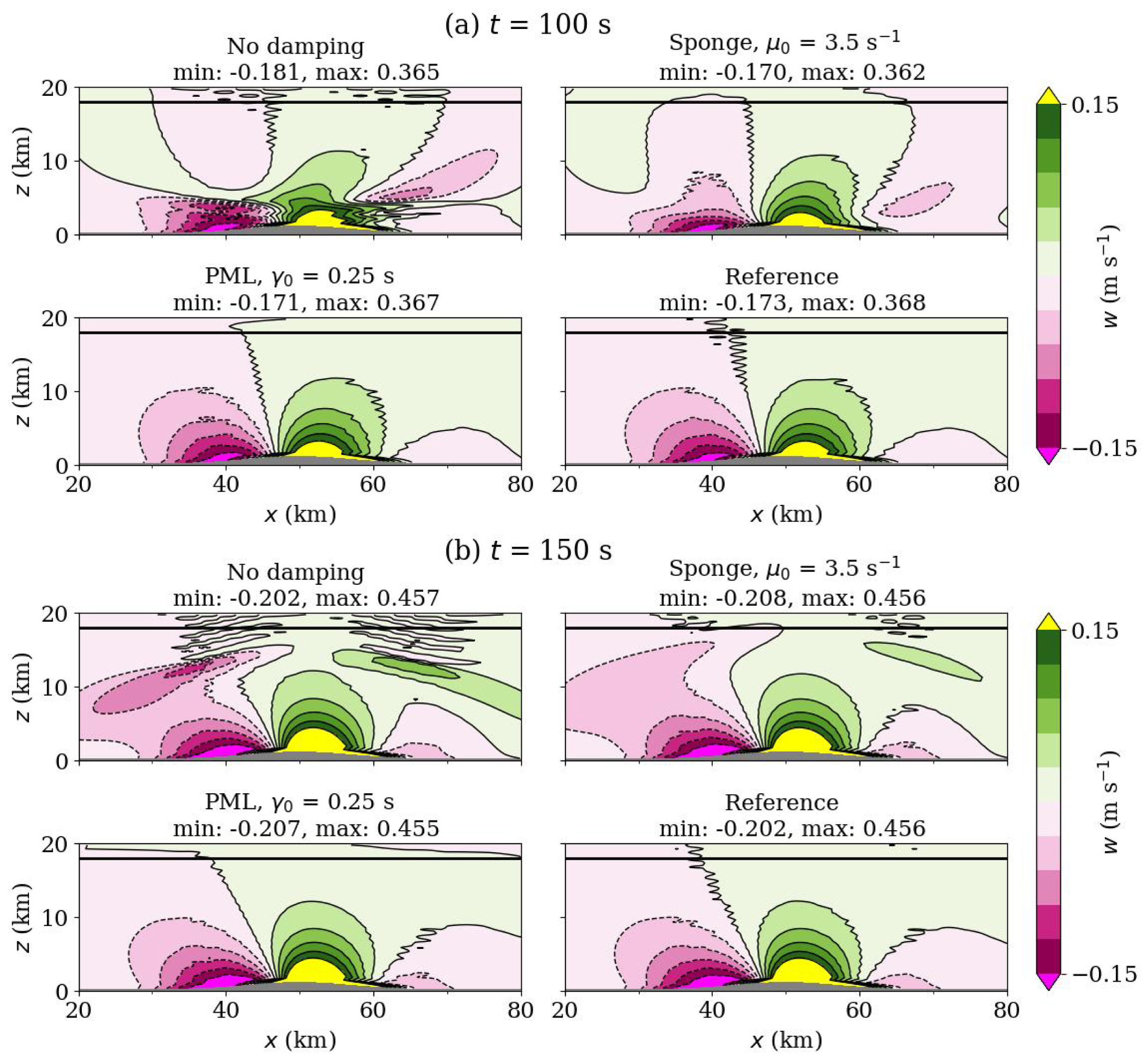}
    \caption{Vertical velocity fields in the nonlinear Boussinesq acoustic test after 100 s (a) and the end time of 150 s (b). The mountain orography is shaded in grey. Solutions are compared from simulations without model top damping, a sponge layer, a stretched PML, and a reference solution with a higher model top.}
    \label{fig:bous_acoustic_w}
\end{figure}

Figure \ref{fig:bous_acoustic_w} shows vertical velocity fields at $t=100 ~\text{s}$ and the final time of $t=150 ~\text{s}$ in the acoustic test. Without damping, there is clear evidence of acoustic wave reflection: at $t=100 ~\text{s}$ the wave has reflected off the model top and interferes with the flow near the mountain, whilst at $t=150 ~\text{s}$ the wave has reflected off the bottom surface and is travelling upwards again. The use of a sponge layer of $\mu_0=3.5 ~\text{s}^{-1}$, which is the choice of damping timescale with the lowest vertical velocity error in this test, damps the acoustic wave to some degree, but there remain spurious features in the solution from wave reflection. With the PML, the acoustic waves are successfully damped, and the flow appears very close to the reference solution. \par

\begin{figure}
    \centering
    \includegraphics[width=0.7\linewidth]{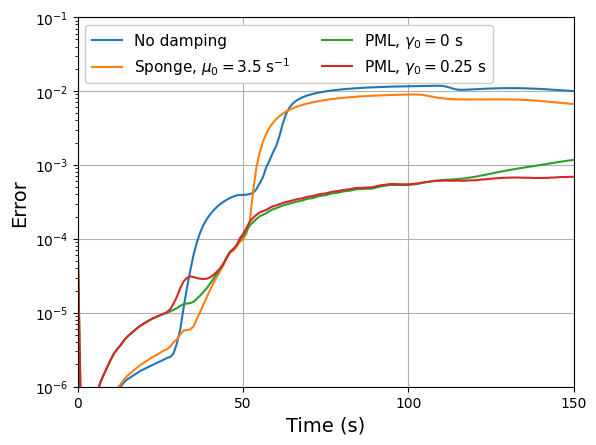}
    \caption{$L_2$ vertical velocity error in the undamped domain, $z \in [0,18] ~\text{km}$, for the nonlinear Boussinesq acoustic test. Comparisons are made of simulations with no model top damping, a sponge layer, an unstretched PML with $\gamma_0=0 ~\text{s}$, and a stretched PML with $\gamma_0 = 0.25 ~\text{s}$.}
    \label{fig:bous_acoustic_error}
\end{figure}

A time series of the $L_2$ vertical velocity error is shown in Figure \ref{fig:bous_acoustic_error}. There is an order of magnitude accuracy improvement with a PML compared to the sponge and undamped solutions. Like with the linear Boussinesq test, a nonzero vertical grid stretching factor reduces the error, although there is minimal difference between $\gamma_0=0 ~\text{s}$ and $\gamma_0=0.25 ~\text{s}$ until $t > 110 ~\text{s}$. Grid stretching factors of $\gamma_0 \in \{0.25, 0.5, 1\} ~\text{s}$ lead to similar solutions. \par

\FloatBarrier

\subsubsection{Orographic gravity wave test results}
We now move to the longer simulation that tests the PML's ability to damp orographic gravity waves. An approximate reference solution is computed using a 50 km model top and a 30 km sponge layer, using a strength of $\mu_0 = 3.5 ~\text{s}^{-1}$. The increased sponge thickness leads to smaller gradients in the damping profile (\ref{eq:damp_coeff}) compared to a thickness of $\delta = 2 ~\text{km}$. During the simulation length, the orographic gravity wave will not travel far enough to reflect off the high model top, and the large sponge should sufficiently damp the acoustic waves. \par
\begin{figure}
    \centering
    \includegraphics[width=0.8\linewidth]{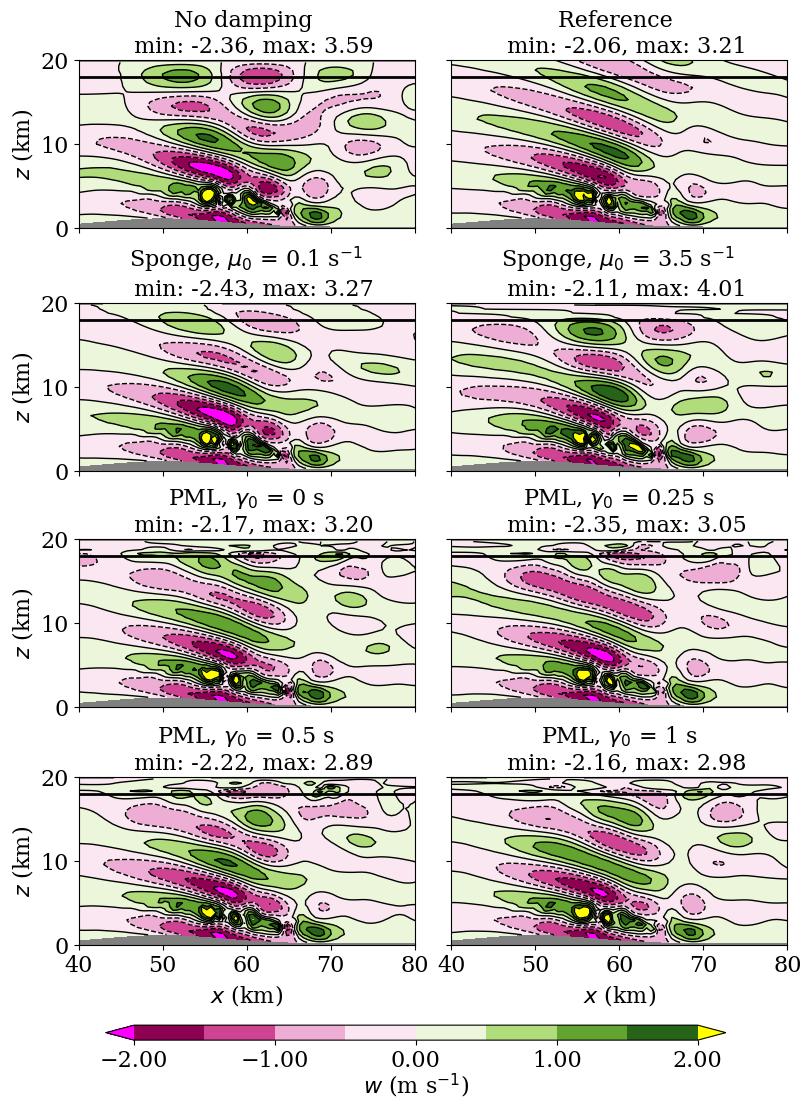}
    \caption{Vertical velocities in the nonlinear Boussinesq orographic gravity wave test after 10,000 s. The black line at $z=18 ~\text{km}$ shows the separation between the damping layer of $z \in [18,20] ~\text{km}$ (for the PML and sponge solutions) and the undamped domain. Field extrema are computed in the undamped domain. Solutions are shown for a reference solution and simulations without model top damping, two strengths of sponge, and PMLs with different real grid stretching factors, $\gamma_0$.}
    \label{fig:bous_gravity_final}
\end{figure}

Figure \ref{fig:bous_gravity_final} compares vertical velocity fields at the end time of 10,000 s. The undamped solution contains a `checkerboard' standing wave pattern due to reflection at the model top, which degrades the representation of the orographic gravity wave throughout the domain. When using a sponge layer, there is a reduction of this checkerboarding, but it is still present. A weaker sponge coefficient is needed compared to the acoustic test, as $\mu_0 = 3.5 ~\text{s}^{-1}$ leads to strong reflections at the onset of the damping layer at $z=18 ~\text{km}$. This highlights a drawback of sponge layers, in that they are difficult to tune to damp both acoustic and gravity waves; this means that sponge layers are typically only designed to damp gravity waves. We identify $\mu_0 = 0.1 ~\text{s}^{-1}$ as a suitable choice for damping the orographic gravity wave in this test. \par
Simulations with a PML in the orographic gravity wave test are visually free of the checkerboard reflection pattern, with a clearer representation of the orographic gravity wave signal for $z \leq 18 ~\text{km}$. However, the use of different grid stretching factors $\gamma_0 \in \{0, 0.25, 0.5, 1 \}$ s modifies the representation of the orographic gravity wave; each choice of $\gamma_0$ replicates its general structure, but there are notable discrepancies relative to the reference solution. Hence, the currently formulated PML successfully avoids the checkerboarding reflection but still needs improvement for a more accurate representation of orographic gravity waves. This will likely require additional terms in the PML equations pertaining to advection or buoyancy. \par

\begin{figure}
    \centering
    \includegraphics[width=\linewidth]{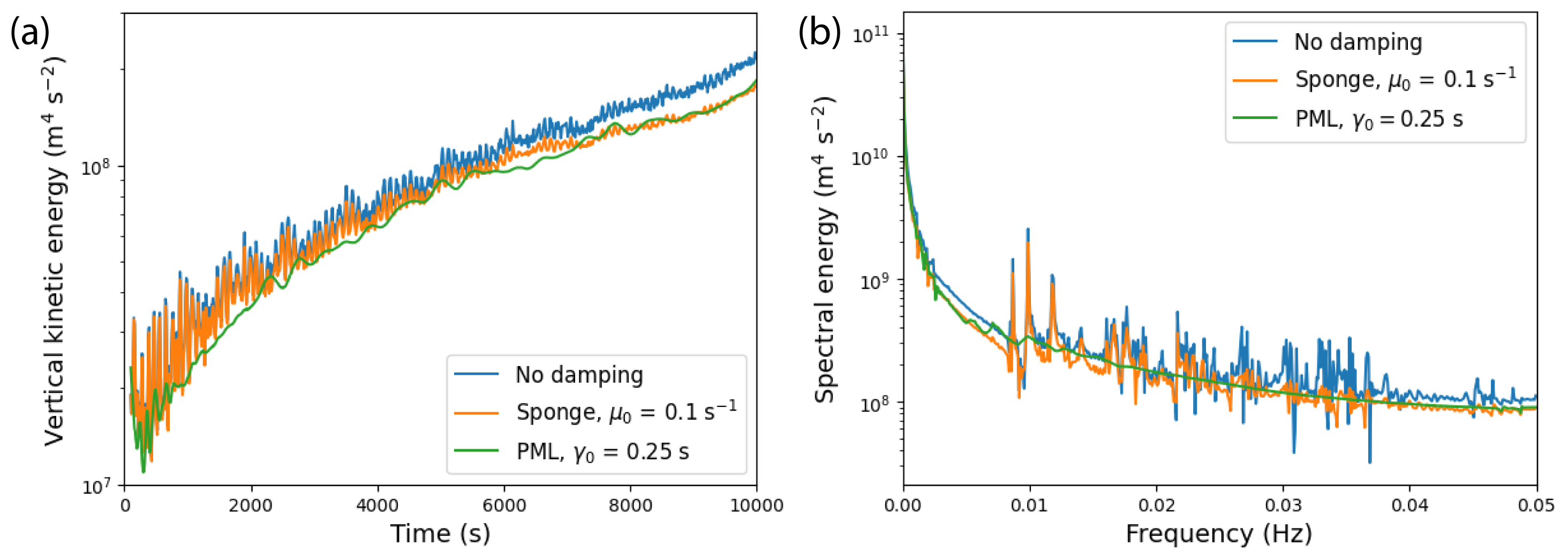}
    \caption{Investigating the integrated vertical kinetic energy of $E_{\text{K,V}} = \intOmega{0.5 w^2}$ in the nonlinear Boussinesq orographic gravity wave test, for simulations with no model top damping, a sponge of $\mu_0=0.1 ~\text{s}^{-1}$, and a PML of $\gamma_0=0.25$ s. (a) A time series of $E_{\text{K,V}}$. (b) The energy spectrum of $|\hat{E}_{\text{K,V}}|$, where the hat denotes the temporal Fourier transform.}
    \label{fig:bous_ke}
\end{figure}

This test still excites an acoustic wave, albeit one that is not properly resolved, due to the implicit timestepping method and a long timestep size of acoustic Courant number $\textrm{Cr}_{\text{ac}}=7$. To investigate the impact of the acoustic waves, Figure \ref{fig:bous_ke}a plots time series of the integrated vertical kinetic energy, $E_{\text{K,V}} = \intOmega{0.5w^2}$. In the undamped and sponge solutions, there are high frequency oscillations in the vertical kinetic energy corresponding to the acoustic wave bouncing between the model top and bottom surface. Such oscillations in the vertical kinetic energy have also been observed for energy-conserving discretisations of the compressible Euler equations on the Lorenz grid \parencite{Lee21,LeePalha21}. The use of the PML leads to a smoother time evolution of $E_{\text{K,V}}$, indicating an absence of this trapped acoustic wave. Hence, the PML can damp the vertically propagating waves even when they are not resolved by the timestepping method. A plot of spectral energy (Figure \ref{fig:bous_ke}b), computed with the temporal Fourier transform of $E_{\text{K,V}}$, shows that there are high frequency components in the undamped and sponge solutions that are removed with the PML. \par

We lastly examine the impact of the PML on the energy in the nonlinear Boussinesq equations. Whilst we have not identified a conserved energy for the nonlinear system, we nevertheless use the same PML power term as for the linear system (\ref{eq:P_PML}) as an approximate measure of energy removal from the PML. We note that for all cases, irrespective of damping choice, there is energy dissipation in the nonlinear equations from upwinding of the advective terms (\ref{eq:nl_bous_mom_weak}). \par
\begin{figure}[htpb]
    \centering
    \includegraphics[width=0.75\linewidth]{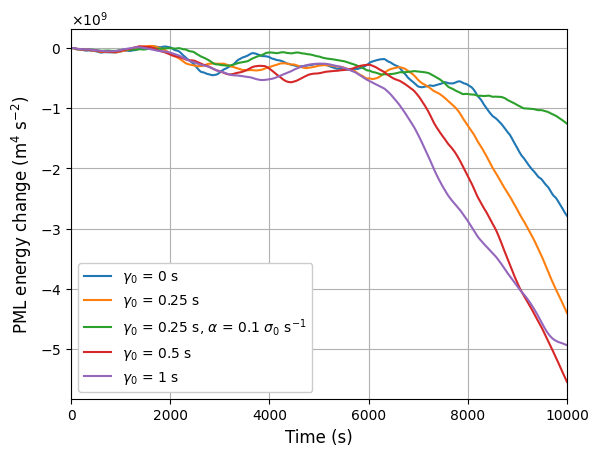}
    \caption{The energy change from the PML in the nonlinear Boussinesq orographic gravity wave test, for different real grid stretching factor, $\gamma_0$. For $t >$ 2000 s, the PML leads to a net energy removal in all cases. Increasing the CFS parameter $\alpha$ reduces the energy removed, as shown for $\gamma_0 = 0.25 ~\text{s}$.}
    \label{fig:bous_gravity_wave_P_PML}
\end{figure}

Figure \ref{fig:bous_gravity_wave_P_PML} examines the net energy removal from the PML by plotting the time-integrated PML power, approximated as
\begin{equation}
    \Delta E_{\text{PML}}(t_n) = \int_{t=0}^{t_n} ~P_{\text{PML}} ~\textrm{d} t \approx \sum_{i=1}^{n} \Delta t ~P_{\text{PML}}(t_i),
\end{equation}

\noindent where $t_i$ are the times used in the discrete solution. Figure \ref{fig:bous_gravity_wave_P_PML} shows that a larger grid stretching factor $\gamma_0$ leads to more energy removal during the simulation; the difference is clearest at the later times, when the orographic gravity wave has reached the PML. Increasing the complex frequency shift factor, $\alpha$, reduces the total energy removed.

\FloatBarrier

\section{Conclusions and discussion}
\label{sec:conclusions}

This paper introduced a novel application of the perfectly matched layer for damping vertically propagating waves in atmospheric models. In the continuous equations, the PML acts to damp all waves that enter the damping region, irrespective of angle of incidence and wavelength, and avoids reflection at the onset of the damping layer. This study derived PML equations for the compressible Boussinesq equations, which are a valuable model of vertical dynamics in the atmosphere and support acoustic and gravity waves. \par
Numerical tests of the PML equations were performed with a compatible finite element model. We set the PML to cover four vertical levels of next-to-lowest-order elements, so although auxiliary PML variables need to be introduced, the increase in computational cost is somewhat mitigated by the thin damping layer. The proposed PML damps motions created by the pressure gradient and divergence terms in the Boussinesq system. Accordingly, this requires velocity and pressure PML variables but not a buoyancy PML variable. \par
Tests with the linear Boussinesq equations showed that the PML is more effective at damping acoustic waves than a standard vertical velocity sponge layer. We then provided a novel use of linear PML theory to the nonlinear Boussinesq equations, which damps perturbations from a hydrostatically balanced state. Numerical tests in the nonlinear system showed that the PML can avoid model top reflections from both acoustic and orographic gravity waves. This is in contrast to the sponge, which still permitted model top reflections in the acoustic test with a strong coefficient, whilst needing a much weaker coefficient for the orographic gravity wave test. \par
The use of a real grid stretching factor in the PML improved the damping of acoustic waves for both the linear and nonlinear Boussinesq systems. However, different stretch factors led to modified representations of the orographic gravity wave in the nonlinear test, although all choices successfully avoided model top reflection. Hence, an important next step is improving our PML for orographic gravity wave applications. This could be achieved by including additional terms in the PML, such as a measure of the buoyancy oscillation frequency or a linearisation of the advective terms, possibly about a nonzero mean flow. Additionally, improvements may also be obtained through different selections of PML parameters. \par
There are a number of future research topics following this study. First is the development of stability proofs for our Boussinesq PML formulation, using error estimates within the Laplace domain. Second, we aim to derive a PML for the compressible Euler equations in the vertical slice domain. The Euler equations are more complex than the Boussinesq equations, due to a nonlinear pressure gradient term and a nonconstant acoustic wave speed. For the compressible Euler PML, we would like to incorporate a linearisation of internal gravity wave propagation. Third, more sophisticated timestepping methods could be used with the PML equations. Specifically, the PML could be incorporated into the semi-implicit quasi-Newton method that is often used for the compressible Euler equations with compatible finite element models, e.g. \textcite{bendall2020compatible}. Last, we eventually aim to implement the PML in a fully three-dimensional model on the sphere, thus testing the potential for a PML to replace Rayleigh damping or Laplacian sponge layers in dynamical cores. \par

\section*{Code Availability}
\sloppy
The Gusto code base is available at the GitHub repository \mbox{\url{https://github.com/firedrakeproject/gusto}}. Python scripts used to run the simulations and create the figures are available at the GitHub repository \mbox{\url{https://github.com/ta440/vertical_slice_PML}}. Functions from the tomplot Python code library were used for the figures, and these can be found at the GitHub repository \mbox{\url{https://github.com/tommbendall/tomplot}}.

\section*{Acknowledgments}
TCA thanks Christiane Jablonowski for the discussions about Rayleigh damping sponge layers, which led him down the path of exploring other options like the PML.

\clearpage

\sloppy

\printbibliography

\end{document}